\documentstyle[11pt]{article}

\author{Charles H. Walter\thanks{Supported in part by NSA research grant
MDA904-92-H-3009.}\\
Laboratoire de Math\'ematiques\\
Universit\'e de Nice\\
F-06108 Nice Cedex 02 FRANCE\\
walter@math.unice.fr}
\title{Components of the Stack of Torsion-Free Sheaves of Rank 2 on Ruled
Surfaces
}
\date{ }

\newtheorem{theorem}[subsection]{Theorem}
\newtheorem{lemma}[subsection]{Lemma}
\newtheorem{proposition}[subsection]{Proposition}

\newtheorem{rk}[subsection]{\sl Remark}
\newenvironment{remark}{\begin{rk}\rm}{\end{rk}}

\input amssym.def
\newsymbol\twoheadrightarrow 1310
\newsymbol\boxtimes 1202
\newsymbol\subsetneqq 2324
\newsymbol\ncong 231D

\def\qed{\hfill $\Box$}
\long\def\TeXButton#1#2{#2}
\def\proof{\paragraph{Proof. }}

\catcode`\@=11

\def\@begintheorem#1#2{\sl \trivlist
   \item[\hskip \labelsep{\bf #1\ #2\thmcounterend}]}
\def\@opargbegintheorem#1#2#3{\sl \trivlist
      \item[\hskip \labelsep{\bf #1\ #2\ (#3)\thmcounterend}]}
\def\thmcounterend{.}

\def\section{\@startsection{section}{1}{\z@}{-3.25ex plus
 -1ex minus -.2ex}{1.5ex plus .2ex}{\large\bf}}
\def\subsection{\@startsection
 {subsection}{2}{\z@}{3.25ex plus 1ex minus .2ex}{-0.5em}{\normalsize\sl}}
\def\subsubsection{\@startsection
 {subsubsection}{3}{\z@}{3.25ex plus 1ex minus .2ex}{-0.5em}{\normalsize\sl}}
\def\paragraph{\@startsection
 {paragraph}{3}{\z@}{2ex plus 0.6ex minus .2ex}{-0.5em}{\normalsize\sl}}
\def\subparagraph{\@startsection
 {subparagraph}{3}{\parindent}{2ex plus 0.6ex minus
.2ex}{-1pt}{\normalsize\sl}}

\def\RIfM@{\relax\ifmmode}

\newif\iffirstchoice@
\firstchoice@true
\def\textfonti{\the\textfont\@ne}
\def\textfontii{\the\textfont\tw@}
\def\text{\RIfM@\expandafter\text@\else\expandafter\text@@\fi}
\def\text@@#1{\leavevmode\hbox{#1}}
\def\text@#1{\mathchoice
 {\hbox{\everymath{\displaystyle}\def\textfonti{\the\textfont\@ne}
\def\textfontii{\the\textfont\tw@}\textdef@@ T#1}}
 {\hbox{\firstchoice@false
  \everymath{\textstyle}\def\textfonti{\the\textfont\@ne}
\def\textfontii{\the\textfont\tw@}\textdef@@ T#1}}
 {\hbox{\firstchoice@false
  \everymath{\scriptstyle}\def\textfonti{\the\scriptfont\@ne}
\def\textfontii{\the\scriptfont\tw@}\textdef@@ S\rm#1}}
 {\hbox{\firstchoice@false
  \everymath{\scriptscriptstyle}\def\textfonti
  {\the\scriptscriptfont\@ne}
\def\textfontii{\the\scriptscriptfont\tw@}\textdef@@ s\rm#1}}}
\def\textdef@@#1{\textdef@#1\rm\textdef@#1\bf\textdef@#1\sl\textdef@#1\it}
\def\DN@{\def\next@}
\def\eat@#1{}
\def\textdef@#1#2{%
 \DN@{\csname\expandafter\eat@\string#2fam\endcsname}%
 \if S#1\edef#2{\the\scriptfont\next@\relax}%
 \else\if s#1\edef#2{\the\scriptscriptfont\next@\relax}%
 \else\edef#2{\the\textfont\next@\relax}\fi\fi}

\catcode`\@=12

\begin{document}

\maketitle
\begin{abstract}
\noindent Let $S$ be a ruled surface without sections of negative
self-intersection. We classify the irreducible components of the moduli
stack of torsion-free sheaves of rank $2$ sheaves on $S$. We also classify
the irreducible components of the Brill-Noether loci in ${\rm Hilb}^N({\bf P}%
^1\times {\bf P}^1)$ given by $W_N^0(D)=\{[X]\mid h^1({\cal I}_X(D))\geq 1\}$
for $D$ an effective divisor class. Our methods are also applicable to ${\bf %
P}^2$ giving new proofs of theorems of Str\o mme (slightly extended) and
Coppo. \bigskip\
\end{abstract}

Let $\pi {:}~S={\bf P}({\cal A})\rightarrow C$ be a ruled surface with
tautological line bundle ${\cal O}(1):={\cal O}_{{\bf P}({\cal A})}(1)$. The
current classification of isomorphism classes of rank $2$ vector bundles $%
{\cal E}$ on $S$ (\cite{BS} \cite{B} \cite{HS} \cite{Ho}) proceeds by
stratifying the moduli functor (or stack) and then classifying the sheaves
in each stratum independently. The numerical data used to distinguish the
strata are usually (i) the splitting type ${\cal O}_{{\bf P}^1}(a)\oplus
{\cal O}_{{\bf P}^1}(b)$ of the generic fiber of $\pi $ (with $a\geq b$),
and (ii) the degree of the locally free sheaf $\pi _{*}({\cal E}(-a))$ on $C$%
. On each stratum, ${\cal U}:=\pi ^{*}(\pi _{*}({\cal E}(-a)))(a)$ is
naturally a subsheaf of ${\cal E}$, and the possible quotient sheaves ${\cal %
E}/{\cal U}$ and extension classes ${\rm Ext}^1({\cal E}/{\cal U},{\cal U})$
have been classified.

To the author's knowledge, rank $2$ torsion-free sheaves on $S$ have not
been given a similar classification, but one could clearly adapt the ideas
used for vector bundles.

What this approach has usually not described is the relationship between the
strata particularly for the strata parametrizing only unstable sheaves. In
this paper we give a first result along these lines by describing which
strata are generic, i.e.\ which are open in the (reduced) moduli stack. Thus
we are really classifying the irreducible components of the moduli stack of
rank $2$ torsion-free sheaves on $S$. We use a method developed by Str\o mme
\cite{St} for rank $2$ vector bundles on ${\bf P}^2$ modified by deformation
theory techniques which originate in \cite{DLP}.

We will divide our irreducible components into two types. The first type we
call prioritary because the general member of a component of this type is a
prioritary sheaf in the sense that we used in \cite{W}. That is, if for each
$p\in C$ we write $f_p:=\pi ^{-1}(p)$ for the corresponding fiber, then a
coherent sheaf ${\cal E}$ on $S$ is {\em prioritary} if it is torsion-free
and satisfies ${\rm Ext}^2({\cal E},{\cal E}(-f_p))$ for all $p$. We showed
in \cite{W} Lemma 7, that if one polarizes $S$ by an ample divisor $H$ such
that $H\cdot (K_S+f_p)<0$, then $H$-semistable sheaves are prioritary. Thus
the prioritary components should be viewed as playing a role one might
otherwise assign to semistable components. But the condition of priority is
simpler to use than semistability because it does not depend on the choice
of a polarization, and moreover the existence problem has a simpler solution
(particularly in higher rank).

The second type of components are nonprioritary ones.

Our main result is the following. We use the convention that if $D\in {\rm NS%
}(S)$, then ${\cal O}_S(D)$ is the line bundle corresponding to the generic
point of the corresponding component ${\rm Pic}^D(S)$ of the Picard scheme.
This is well-defined on all surfaces for which numerical and algebraic
equivalence coincide, including all of ours.

\begin{theorem}
\label{ruled}Let $S$ be a ruled surface without curves of negative
self-intersection, and let $f$ be the numerical class of a fiber of $S$. Let
$c_1\in {\rm NS}(S)$ and $c_2\in {\bf Z}$. The irreducible components of the
stack ${\rm TF}_S(2,c_1,c_2)$ of torsion-free sheaves on $S$ of rank $2$ and
Chern classes $c_1$ and $c_2$ are the following:

(i)\quad A unique prioritary component if $c_1f$ is even and $c_2\geq \frac
14c_1^2$, or if $c_1f$ is odd and $c_2\in {\bf Z}$. This component is
generically smooth of dimension $-\chi ({\cal E},{\cal E})$, and the general
sheaf in the component is locally free.

(ii)\quad For every pair $(D,n)\in {\rm NS}(S)\times {\bf Z}$ such that $%
Df\leq -1+\frac 12c_1f$ and $0\leq n\leq c_2+D(D-c_1)\leq \chi ({\cal O}%
_S(-c_1))+D(2D-2c_1-K_S)$ a unique nonprioritary component whose general
member is a general extension%
$$
0\rightarrow {\cal I}_{Z_1}(c_1-D)\rightarrow {\cal E}\rightarrow {\cal I}%
_{Z_2}(D)\rightarrow 0
$$
where $Z_1$ (resp.\ $Z_2)$ is a general set of $n$ (resp.\ $n^{\prime
}:=c_2+D(D-c_1)-n$) points on $S$. These components have dimensions $-\chi (%
{\cal E},{\cal E})+\chi ({\cal O}_S(-c_1))+D(D-c_1-K_S)-c_2$ but have
generic embedding codimension $n^{\prime }+h^1({\cal O}_S(2D-c_1)$. The
general sheaf in the component is locally free except at $Z_1$.
\end{theorem}

For ${\bf P}^2$ the components of ${\rm TF}_{{\bf P}^2}(2,c_1,c_2)$
containing locally free sheaves were classified by Str\o mme using a similar
method (\cite{St} Theorem 3.9). We wish to add to his classification the
components of ${\rm TF}_{{\bf P}^2}(2,c_1,c_2)$ whose general member is not
locally free. We recall from \cite{HL} that a prioritary sheaf ${\cal E}$ on
${\bf P}^2$ is one that is torsion-free and satisfies ${\rm Ext}^2({\cal E},%
{\cal E}(-1))=0$.

\begin{theorem}
\label{P2}Let $S$ be ${\bf P}^2$ and let $f\in {\rm NS}(S)$ be the class of
a line. Let $(c_1,c_2)\in {\bf Z}^2$. Then the irreducible components of $%
{\rm TF}_{{\bf P}^2}(2,c_1,c_2)$ have the same classification as in Theorem
\ref{ruled} except that the prioritary component exists if and only if $%
c_2\geq \frac 14c_1^2-\frac 14$.
\end{theorem}

The uniqueness of the prioritary components was proven for ruled surfaces
(resp.\ for ${\bf P}^2$) in \cite{W} (resp.\ \cite{HL}) although of course
there were many earlier results by many authors concerning semistable
components on ${\bf P}^2$ and on various ruled surfaces.

The classification of the irreducible components of the stacks of
torsion-free sheaves has an interesting application to Brill-Noether
problems. Let $S$ be a smooth projective algebraic surface, $E$ an effective
divisor class on $S$, and $N$ a positive integer such that $N\leq h^0({\cal O%
}_S(E))$. For simplicity we will assume that $H^1({\cal O}_S)=H^1({\cal O}%
_S(E))=0$. We consider the Brill-Noether loci in ${\rm Hilb}^NS$ defined as%
$$
W_N^i(E)=\{[X]\in {\rm Hilb}^NS\mid h^1({\cal I}_X(E))\geq i+1\}.
$$
Thus $W_N^i(E)$ parametrizes those $0$-schemes of length $N$ which impose at
least $i+1$ redundant conditions on divisors in $|E|$. What we wish to
consider is:

\paragraph{The Brill-Noether Problem.}

Classify the irreducible components of the $W_N^i(E)$ and compute their
dimensions.\medskip\

It is known from general principles that each component has codimension at
most $(\chi +i+1)(i+1)$ in ${\rm Hilb}^NS$ where $\chi =h^0({\cal O}%
_S(E))-N\geq 0$, but there can be many components of various smaller
codimensions.

The Brill-Noether problem is related to the problem of classifying
irreducible components of the stack of torsion-free sheaves on $S$ as
follows. By an elementary argument (cf.\ \cite{C} p.\ 732) the general $X$
in any component of $W_N^i(E)$ has $h^1({\cal I}_X(E))=i+1$. One then uses
Serre duality $H^1({\cal I}_X(E))^{*}\cong {\rm Ext}^1({\cal I}_X(E-K_S),%
{\cal O}_S)$ to get an extension%
$$
0\rightarrow {\cal O}_S^{\oplus i+1}\rightarrow {\cal E\rightarrow I}%
_X(E-K_S)\rightarrow 0.
$$
So we get a Serre correspondence between $X\in W_N^i(E)-W_N^{i-1}(E)$ and
pairs $({\cal E},V)$ where ${\cal E}$ is torsion-free of rank $i+2$ with $%
c_1({\cal E})=E-K_S$, $c_2({\cal E})=N$, $H^1({\cal E}(K_S))=H^2({\cal E}%
(K_S))=0$, and for which there exists an $(i+1)$-dimensional subspace $%
V\subset H^0({\cal E})$ such that the natural map $V\otimes {\cal O}%
_S\rightarrow {\cal E}$ is injective with torsion-free quotient. Note that
these properties are all open conditions on ${\cal E}$ within the stack of
torsion-free sheaves on $S$. So Theorems \ref{ruled} and \ref{P2} yields a
classification of the irreducible components of the $W_N^0(E)$ for ${\bf P}%
^1\times {\bf P}^1$ and ${\bf P}^2$. This classification has been previously
obtained by Coppo for ${\bf P}^2$ by a different method (\cite{C}
Th\'eor\`eme 3.2.1) but seems new for ${\bf P}^1\times {\bf P}^1$.

\begin{theorem}
\label{BN}Let $S$ be ${\bf P}^2$ (resp.\ ${\bf P}^1\times {\bf P}^1$), let $%
E $ be an effective divisor of degree $e$ (resp.\ of bidegree $(e_1,e_2)$),
and let $N$ be an integer such that $0<N\leq \chi ({\cal O}_S(E))$. Then the
irreducible components of the Brill-Noether locus $W_N^0(E)$ are the
following:

(i)\quad For every pair $(D,n)\in {\rm NS}(S)\times {\bf Z}$ such that $D$
is an effective and irreducible divisor class of degree $d$ on ${\bf P}^2$
(resp.~of bidegree $(d_1,d_2)$ on ${\bf P}^1\times {\bf P}^1$) such that $%
d\leq \frac 12(e+1)$\ (resp.\ $d_2\leq \frac 12e_2$), $D(E-D)\leq \chi (%
{\cal O}_S(E))-N$, $n\geq 0$, and $0\leq N-D(E-D-K_S)-n\leq \chi ({\cal O}%
_S(D+K_S))$, there exists a unique irreducible component of codimension $%
D(E-D)+1$ in ${\rm Hilb}^N(S)$ whose general member is the union of $n$
general points of $S$ and $N-n$ points on a curve in $\left| D\right| $.

(ii)\quad If $S$ is ${\bf P}^2$ (resp.\ if $S$ is ${\bf P}^1\times {\bf P}^1$
and $e_2$ is even, resp.\ if $S$ is ${\bf P}^1\times {\bf P}^1$ and $e_2$ is
odd), then there exists one additional component of codimension $\chi ({\cal %
O}_S(E))$$-N+1$ in ${\rm Hilb}^N(S)$ if $N\geq \frac 14(e+2)(e+4)$ (resp.~$%
N\geq \frac 12(e_1+2)(e_2+2)$, resp.\ $N\geq \frac 12(e_1+2)(e_2+1)+1$). If $%
S={\bf P}^1\times {\bf P}^1$ and $(e_1,e_2,N)=(e_1,1,e_1+2)$ there is also
one additional component of codimension $\chi ({\cal O}_S(E))$$-N+1$.
\end{theorem}

In part (i) the $N-n$ points on the curve $C\in |D|$ have the property that
their union is a divisor on $C$ belonging to a linear system of the form $%
\left| \Gamma +E{\mid }_C-K_C\right| $ with $\Gamma $ an effective divisor
satisfying $h^0({\cal O}_C(\Gamma ))=1$. The necessary condition $0\leq \deg
(\Gamma )\leq g(C)$ is exactly the condition $0\leq N-D(E-D-K_S)-n\leq \chi (%
{\cal O}_S(D+K_S))$.

The main tool which we use to obtain our results is interesting in its own
right. We use the notation $\chi ({\cal F},{\cal G})=\sum (-1)^i\dim {\rm Ext%
}^i({\cal F},{\cal G})$.

\begin{proposition}
\label{unobst}Let $S$ be a projective surface, and ${\cal E}$ a coherent
sheaf on $S$ with a filtration $0=F_0({\cal E})\subset F_1({\cal E})\subset
\cdots \subset F_r({\cal E})={\cal E}$. Suppose that the graded pieces ${\rm %
gr}_i({\cal E}):=F_i({\cal E})/F_{i-1}({\cal E})$ satisfy ${\rm Ext}^2({\rm %
gr}_i({\cal E}),{\rm gr}_j({\cal E}))$ for $i\geq j$. Then

(i)\quad the deformations of ${\cal E}$ as a filtered sheaf are unobstructed,

(ii)\quad if ${\cal E}$ is a generic filtered sheaf, then the ${\rm gr}_i(%
{\cal E})$ are generic, and

(iii)\quad if ${\cal E}$ is generic as an unfiltered sheaf, then also $\chi (%
{\rm gr}_i({\cal E}),{\rm gr}_{i+1}({\cal E}))\geq 0$ for $i=1,\ldots ,r-1$.
\end{proposition}

The outline of the paper is as follows. In the first section we review some
necessary facts about algebraic stacks and their dimensions. In the second
section we prove our technical tool Proposition \ref{unobst} and describe
some situations where it applies. It the third section we classify the
prioritary components of the ${\rm TF}_S(2,c_1,c_2)$ and the $W_N^0(E)$. In
the fourth section we classify the nonprioritary components. In the short
final section we complete the proofs of the main theorems.

This paper was written in the context of the group on vector bundles on
surfaces of Europroj. The author would like to thank A.\ Hirschowitz and
M.-A.\ Coppo for some useful conversations.

\section{Algebraic Stacks}

In this paper we use stacks because in that context there exist natural
universal families of coherent (or torsion-free) sheaves. The paper should
be manageable even to the reader unfamiliar with algebraic stacks if he
accepts them as some sort of generalization of schemes where there are
decent moduli for unstable sheaves. For the reader who wishes to learn about
algebraic stacks we suggest \cite{LMB}. Alternative universal families of
coherent sheaves which stay within the category of schemes would be certain
standard open subschemes of Quot schemes. This is the approach taken in \cite
{St}. But the language of algebraic stacks is the natural one for problems
which involve moduli of unstable sheaves.

Stacks differ from schemes is in the way their {\em dimensions} are
calculated. For the general definition of the dimension of an algebraic
stack at one of its points the reader should consult \cite{LMB} \S 5. But
the dimension of the algebraic stack ${\rm Coh}_S$ of coherent sheaves on $S$
(or of any open substack of ${\rm Coh}_S$ such as a ${\rm TF}_S(r,c_1,c_2)$)
at a point corresponding to a sheaf ${\cal E}$ is the dimension of the
Kuranishi formal moduli for deformations of ${\cal E}$ (i.e.\ the fiber of
the obstruction map $({\rm Ext}^1({\cal E},{\cal E}),0)^{\wedge }\rightarrow
({\rm Ext}^2({\cal E},{\cal E}),0)^{\wedge }$) minus the dimension of the
automorphism group of ${\cal E}$. Thus if we write $e_i=\dim {\rm Ext}^i(%
{\cal E},{\cal E})$, then $-e_0+e_1-e_2\leq \dim {}_{[{\cal E}]}{\rm Coh}%
_S\leq -e_0+e_1$. If $S$ is a surface, this means
\begin{equation}
\label{dimcoh}-\chi ({\cal E},{\cal E})\leq \dim {}_{[{\cal E}]}{\rm Coh}%
_S\leq -\chi ({\cal E},{\cal E})+e_2.
\end{equation}
If ${\cal E}$ is a stable sheaf, then $\dim {}_{[{\cal E}]}{\rm Coh}_S$ is
one less than the dimension of the moduli scheme at $\left[ {\cal E}\right] $
because ${\cal E}$ has a one-dimensional family of automorphisms, the
homotheties.

Generally speaking, the dimension of an algebraic stack are well-behaved. It
is constant on an irreducible component away from its intersection with
other components; the dimension of a locally closed substack is smaller than
the dimension of the stack; etc. But stacks can have negative dimensions.

\section{When is a Filtered Sheaf Generic?}

In this section we prove our main technical tool Proposition \ref{unobst}
and then give two corollaries applying the proposition to birationally ruled
surfaces.

\paragraph{Proof of Proposition \ref{unobst}.}

We begin by recalling some of the deformation theory of \cite{DLP}. We
consider the abelian category of sheaves with filtrations of a fixed length $%
r$:%
$$
0=F_0({\cal E})\subset F_1({\cal E})\subset \cdots \subset F_r({\cal E})=%
{\cal E.}
$$
On this category we can define functors%
$$
\begin{array}{rcl}
{\rm Hom}_{-}({\cal E},{\cal F}) & = & \{\phi \in
{\rm Hom}({\cal E},{\cal F})\mid \phi (F_i({\cal E}))\subseteq F_i({\cal F})
\text{ for all }i\}, \\ {\rm Hom}_{neg}({\cal E},{\cal F}) & = & \{\phi \in
{\rm Hom}({\cal E},{\cal F})\mid \phi (F_i({\cal E}))\subseteq F_{i-1}({\cal %
F})\text{ for all }i\}.
\end{array}
$$
These have right-derived functors denoted ${\rm Ext}_{-}^p$ and ${\rm Ext}%
_{neg}^p$ which may be computed by the spectral sequences (\cite{DLP}
Proposition 1.3)%
\begin{eqnarray}
E_1^{pq} \  = & \left\{
\begin{array}{ll}
\prod_i{\rm Ext}^{p+q}({\rm gr}_i({\cal E}),{\rm gr}_{i-p}({\cal E})) &
\text{if }p\geq 0 \\ 0 & \text{if }p\le -1
\end{array}
\right\} & \Rightarrow  \  {\rm Ext}_{-}^{p+q}
({\cal E},{\cal E}), \label{specf} \\
E_1^{pq} \  = & \left\{
\begin{array}{ll}
\prod_i{\rm Ext}^{p+q}({\rm gr}_i({\cal E}),{\rm gr}_{i-p}({\cal E})) &
\text{if }p\geq 1 \\ 0 & \text{if }p\leq 0
\end{array}
\right\} & \Rightarrow \  {\rm Ext}_{neg}^{p+q}({\cal E},{\cal E}).
\label{specminus}
\end{eqnarray}
There is also a long exact sequence
\begin{equation}
\label{extneg}\cdots \rightarrow {\rm Ext}_{neg}^p({\cal E},{\cal E}%
)\rightarrow {\rm Ext}_{-}^p({\cal E},{\cal E})\rightarrow \prod_i{\rm Ext}%
^p({\rm gr}_i({\cal E}),{\rm gr}_i({\cal E}))\rightarrow {\rm Ext}%
_{neg}^{p+1}({\cal E},{\cal E})\rightarrow \cdots .
\end{equation}

(i) The tangent space for the deformations of ${\cal E}$ as a filtered sheaf
is ${\rm Ext}_{-}^1({\cal E},{\cal E})$ and the obstruction space is ${\rm %
Ext}_{-}^2({\cal E},{\cal E})$. The latter vanishes because of the spectral
sequence (\ref{specf}).

(ii) From (\ref{specminus}) and (\ref{extneg}) we see that the map ${\rm Ext}%
_{-}^1({\cal E},{\cal E})\rightarrow \prod_i{\rm Ext}^1({\rm gr}_i({\cal E}),%
{\rm gr}_i({\cal E}))$ is surjective. Thus any first-order infinitesimal
deformation of the ${\rm gr}_i({\cal E})$ can be induced from a first-order
infinitesimal deformation of ${\cal E}$ as a filtered sheaf. But because of
(i) any first-order infinitesimal deformation of the filtered sheaf ${\cal E}
$ is induced by a noninfinitesimal deformation of ${\cal E}$. So if ${\cal E}
$ is generic, then the ${\rm gr}_i({\cal E})$ must also be generic in their
respective stacks.

(iii) We consider ${\cal E}$ with two filtrations: the original filtration
and its subfiltration obtained by suppressing the term $F_i({\cal E})$. We
write ${\rm Ext}_{-}^p$ (resp.\ ${\rm Ext}_{-,sub}^p$) for the ${\rm Ext}%
_{-}^p$ associated to these two filtrations. We have a long exact sequence
$$
\cdots \rightarrow {\rm Ext}_{-}^p({\cal E},{\cal E})\rightarrow {\rm Ext}%
_{-,sub}^p({\cal E},{\cal E})\rightarrow {\rm Ext}^p({\rm gr}_i({\cal E}),%
{\rm gr}_{i+1}({\cal E}))\rightarrow {\rm Ext}_{-}^{p+1}({\cal E},{\cal E}%
)\rightarrow \cdots .
$$
Also ${\rm Ext}_{-}^2({\cal E},{\cal E})=0$ by (i). So the formal moduli for
the deformations of ${\cal E}$ as a filtered sheaf for the full filtration
is of dimension $\dim {\rm Ext}_{-}^1({\cal E},{\cal E})$. The formal moduli
for the deformations of ${\cal E}$ as a filtered sheaf for the subfiltration
is by general principles of dimension at least%
$$
\dim {\rm Ext}_{-,sub}^1({\cal E},{\cal E})-\dim {\rm Ext}_{-,sub}^2({\cal E}%
,{\cal E})\geq \dim {\rm Ext}_{-}^1({\cal E},{\cal E})-\chi ({\rm gr}_i(%
{\cal E}),{\rm gr}_{i+1}({\cal E})).
$$
So if $\chi $$({\rm gr}_i({\cal E}),{\rm gr}_{i+1}({\cal E}))<0$, then the
natural morphism from the formal moduli for the deformations of ${\cal E}$
with the full filtration to the formal moduli for the deformations of ${\cal %
E}$ with the subfiltration could not be surjective. So there would be finite
deformations of ${\cal E}$ which preserve the subfiltration but not the full
filtration. This would contradict the genericity of ${\cal E}$ as an
unfiltered sheaf.\TeXButton{qed}{\qed \medskip}

There are several situations in which there are filtrations to which
Proposition \ref{unobst} applies. For the first situation, let $S$ be a
smooth projective surface and $H$ an ample divisor on $S$. Recall that the $%
H $-slope of a torsion-free sheaf ${\cal F}$ on $S$ is $\mu _H({\cal F}%
):=(Hc_1({\cal F}))/{\rm rk}({\cal F})$. We write $\mu _{H,\max }({\cal F})$
is the maximum $H$-slope of a nonzero subsheaf of ${\cal F}$, and $\mu
_{H,\min }({\cal F})$ is the minimum slope of a nonzero torsion-free
quotient sheaf of ${\cal F}$.

\begin{lemma}
\label{mu}Let $S$ be a smooth projective surface and $H$ an ample divisor on
$S$ such that $HK_S<0$. Let ${\cal F}$ and ${\cal G}$ be torsion-free
sheaves on $S$ such that $\mu _{H,\max }({\cal F})+HK_S<\mu _{H,\min }({\cal %
G})$. Then ${\rm Ext}^2({\cal F},{\cal G})=0$.
\end{lemma}

\TeXButton{Proof}{\proof}By Serre duality we have ${\rm Ext}^2({\cal F},%
{\cal G})\cong {\rm Hom}({\cal G},{\cal F}(K_S))^{*}$. If there were a
nonzero $\phi \in {\rm Hom}({\cal G},{\cal F}(K_S))$, then we would have
$$
\mu _{H,\max }({\cal F})+HK_S=\mu _{H,\max }({\cal F}(K_S))\geq \mu ({\rm im}%
(\phi ))\geq \mu _{H,\min }({\cal G}),
$$
a contradiction.\TeXButton{qed}{\qed \medskip}

It follows that if $(S,{\cal O}_S(H))$ is a polarized surface such that $%
HK_S<0$, then Proposition \ref{unobst} applies to the Harder-Narasimhan
filtration of a torsion-free sheaf ${\cal E}$ on $S$. It also applies to the
weak Harder-Narasimhan filtration for torsion-free sheaves on ${\bf P}^2$
described in \cite{W2}.

The other situation in which Proposition \ref{unobst} applies is the
relative Harder-Narasimhan filtration of a torsion-free sheaf ${\cal E}$ on
a ruled surface $\pi {:}~S\rightarrow C$. To describe this let $f_\eta $ be
the generic fiber of $\pi $. Write ${\cal E}{\mid }_{f_\eta }\cong
\bigoplus_{i=1}^s{\cal O}_{f_\eta }(e_i)^{n_i}$ with $e_1>e_2>\cdots >e_s$
and the $n_i>0$. There exists a unique filtration $0=F_0({\cal E})\subset
F_1({\cal E})\subset \cdots \subset F_s({\cal E})={\cal E}$ such that the
graded pieces ${\rm gr}_i({\cal E})$ are torsion-free and satisfy ${\rm gr}%
_i({\cal E}){\mid }_{f_\eta }\cong {\cal O}_{f_\eta }(e_i)^{n_i}$. The $F_i(%
{\cal E})$ may be obtained as the inverse image in ${\cal E}$ of the torsion
subsheaf of ${\cal E}/{\cal E}_i$ where ${\cal E}_i$ is the image of the
natural map $\pi ^{*}($$\pi _{*}({\cal E}(-e_i)))(e_i)\rightarrow {\cal E}$.
Proposition \ref{unobst} applies to this relative Harder-Narasimhan
filtration because of

\begin{lemma}
\label{fiber}Let $\pi {:}~S\rightarrow C$ be a ruled surface, and let ${\cal %
E}$ and ${\cal G}$ be torsion-free sheaves on $S$. Suppose that the
restrictions of ${\cal E}$ and ${\cal G}$ to a general fiber $F$ of $\pi $
are of the forms ${\cal E}{\mid }_F\cong \bigoplus_i{\cal O}_F(e_i)$ and $%
{\cal G}{\mid }_F\cong \bigoplus_j{\cal O}_F(g_j)$ with $\max \{e_i\}-2<\min
\{g_j\}.$ Then ${\rm Ext}^2({\cal E},{\cal G})=0$. In particular, if $\max
\{e_i\}-\min \{e_j\}<2$, then ${\cal E}$ is prioritary.
\end{lemma}

\TeXButton{Proof}{\proof}Again Serre duality gives ${\rm Ext}^2({\cal E},%
{\cal G})\cong {\rm Hom}({\cal G},{\cal E}(K_S))^{*}$. If there were a
nonzero $\phi \in {\rm Hom}({\cal G},{\cal E}(K_S))$, then there would be a
nonzero $\phi {\mid }_F\in \bigoplus_{i,j}H^0({\cal O}_F(e_i-2-g_j))$. This
is impossible since $f_i-2-g_j<0$ for all $i$ and $j$.

If $\max \{e_i\}-\min \{e_i\}<2$, then we may set ${\cal G}={\cal E}(-f_p)$
for any fiber $f_p=\pi ^{-1}(p)$ to get ${\rm Ext}^2({\cal E},{\cal E}%
(-f_p))=0$ for all $p\in C$. Thus ${\cal E}$ is prioritary.\TeXButton{qed}
{\qed \medskip}

\section{Prioritary Components}

In this section we prove the necessary lemmas for classifying the principal
components of the ${\rm TF}_S(2,c_1,c_2)$ and $W_N^0(E)$. We use \cite{W}
and \cite{HL} as our basic sources for existence and uniqueness results
because these use our preferred language of prioritary sheaves. But
existence and uniqueness results for only marginally different classes of
sheaves on ${\bf P}^2$ and of rank $2$ sheaves on ruled surfaces had already
been proven in \cite{Ba} \cite{BS} \cite{B} \cite{DLP} \cite{E} \cite{ES}
\cite{HS} \cite{Ho} \cite{Hu1} \cite{Hu2} (and perhaps elsewhere).

\begin{proposition}
\label{discquad}Let $\pi {:}~S\rightarrow C$ be a ruled surface, and let $%
f\in {\rm NS}(S)$ be the numerical class of a fiber of $\pi $. Then ${\rm TF}%
_S(r,c_1,c_2)$ has a unique prioritary component if $r$ divides $c_1f$ and $%
2rc_2\geq (r-1)c_1^2$, or if $r$ does not divide $c_1f$. Otherwise it has no
prioritary components. The prioritary component is smooth of dimension $%
-\chi ({\cal E},{\cal E})$.
\end{proposition}

\TeXButton{Proof}{\proof}The uniqueness and smoothness of the prioritary
component was proven in \cite{W} Proposition 2. Since by definition a
prioritary sheaf ${\cal E}$ satisfies ${\rm Ext}^2({\cal E},{\cal E}%
(-f_p))=0 $ for all $p\in C$, it also satisfies ${\rm Ext}^2({\cal E},{\cal E%
})=0$. So the prioritary component has dimension $-\chi ({\cal E},{\cal E})$
according to (\ref{dimcoh}).

For existence of a prioritary sheaf ${\cal E}$, note that $2rc_2-(r-1)c_1^2$
is invariant under twist as is the residue of $c_1f$ modulo $r$. Then by
replacing ${\cal E}$ by a twist ${\cal E}(n)$ if necessary, we may assume
that $d:=-c_1f$ satisfies $0\leq d<r$. In the proof of \cite{W} Proposition
2 it was shown that a general prioritary sheaf with such a $c_1$ fits into
an exact sequence
\begin{equation}
\label{Beil}0\rightarrow \pi ^{*}({\cal K})\rightarrow {\cal E}\rightarrow
\pi ^{*}({\cal L})\otimes \Omega _{S/C}(1)\rightarrow 0
\end{equation}
where ${\cal K}$ is a vector bundle on $C$ of rank $r-d$ and ${\cal L}$ a
coherent sheaf on $C$ of rank $d$. Let $k=\deg ({\cal K})$ and $l=\deg (%
{\cal L})$. Write $h=c_1({\cal O}(1))$ so that $\{h,f\}$ is a basis of ${\rm %
NS}(S)$. Then ${\cal E}$ has rank $r$ and Chern classes $c_1=(k+l)f-dh$ and $%
c_2=\frac 12d(d-1)h^2-(k+l)d+l$. So to finish the proof of the lemma we need
to show that if $0<d<r$, then there exist prioritary sheaves of the form (%
\ref{Beil}) for all $k$ and $l$, while if $d=0$, then there exist prioritary
sheaves of that form if and only if $(k,l)$ satisfies $l\geq 0$.

If $0<d<r$, then for any $k$ and $l$ and any locally free sheaves ${\cal K}$
(resp.\ ${\cal L})$ on $C$ of rank $r-d$ and degree $k$ (resp.\ rank $d$ and
degree $l$), the sheaf ${\cal F}:=\pi ^{*}({\cal K})\oplus \left[ \pi ^{*}(%
{\cal L})\otimes \Omega _{S/C}(1)\right] $ has splitting type ${\cal O}_{%
{\bf P}^1}^{r-d}\oplus {\cal O}_{{\bf P}^1}(-1)^d$ on all fibers and hence
is prioritary by Lemma \ref{fiber}.

If $d=0$, then ${\cal L}$ has rank $0$. So its degree $l$ must be
nonnegative. Conversely if $k$ is any integer and $l\geq 0$, then an ${\cal E%
}$ as in (\ref{Beil}) can be constructed for any locally free sheaf ${\cal K}
$ of rank $r$ and degree\thinspace $k$ on $C$ as an elementary transform of $%
\pi ^{*}({\cal K})$ along $l$ fibers of $\pi $. Such an ${\cal E}$ is
prioritary by Lemma \ref{fiber} because its restriction to the general fiber
of $\pi $ is trivial. Thus for $d=0$ there exists prioritary sheaf ${\cal E}$
of the form (\ref{Beil}) for and only for those $(k,l)$ satisfying $l\geq 0$%
. This completes the proof of the lemma. \TeXButton{qed}{\qed \medskip}

\begin{proposition}
\label{discp2}{\rm (Hirschowitz-Laszlo)} The stack ${\rm TF}_{{\bf P}%
^2}(r,c_1,c_2)$ has a unique prioritary component if $2rc_2-(r-1)c_1^2\geq
-d(r-d)$ where $c_1\equiv -d\pmod{r}$ and $0\leq d<r$. Otherwise it has no
prioritary components. The prioritary component is smooth of dimension $%
-\chi ({\cal E},{\cal E})$.
\end{proposition}

\TeXButton{Proof}{\proof}Let $c_1=mr-d$. Let $\mu =c_1/r$ be the slope and $%
\Delta =(2rc_2-(r-1)c_1^2)/2r^2$ the discriminant of ${\cal E}$. Then in
\cite{HL} Chap.~I, Propositions 1.3 and 1.5 and Th\'eor\`eme 3.1, it is
shown that ${\rm TF}_{{\bf P}^2}(r,c_1,c_2)$ has a priority component if and
only if the Hilbert polynomial $P(n)=r\left( \frac 12(\mu +n+2)(\mu
+n+1)-\Delta \right) $ is nonpositive for some integer $n$, and that in that
case the prioritary component is unique and smooth. The dimension of such a
component is again $-\chi ({\cal E},{\cal E})$ by (\ref{dimcoh}) because the
prioritary condition implies ${\rm Ext}^2({\cal E},{\cal E})=0$.

We show that the Hilbert polynomial criterion of \cite{HL} is equivalent to
the criterion asserted by the lemma. But $P(n)-P(n-1)=\mu +n+1=m-\frac
dr+n+1 $ is nonnegative if and only if $n\geq -m$. So $\min _{n\in {\bf Z}%
}P(n)=P(-m-1)=r\left( \frac 12(1-\frac dr)(-\frac dr)-\Delta \right) $, and
this is nonpositive if and only if $2r^2\Delta \geq -d(r-d)$.\thinspace
\TeXButton{qed}{\qed \medskip}

We recall the Riemann-Roch formula for a coherent sheaf ${\cal E}$ of rank $%
r $ and Chern classes $c_1$ and $c_2$ on a surfaces $S$:
\begin{equation}
\label{RR}\chi ({\cal E})=r\chi ({\cal O}_S)+\frac 12c_1\left(
c_1-K_S\right) -c_2.
\end{equation}

\begin{lemma}
\label{diminution}Let $\pi {:}~S\rightarrow C$ be a ruled surface or let $S$
be ${\bf P}^2$. Suppose ${\cal E}$ is a prioritary sheaf on $S$ of rank $%
r\geq 2$ such that $H^1({\cal E})=H^2({\cal E})=0$. Let $H$ be a very ample
divisor on $S$.

(i)\quad If ${\cal F}$ is a general prioritary sheaf of rank $r$ and Chern
classes $c_1=c_1({\cal E})$ and $c_2\geq c_2({\cal E})$ such that $\chi (%
{\cal F})\geq 0$, then $H^1({\cal F})=H^2({\cal F})=0$.

(ii)\quad If in addition $H^1({\cal E}(H))=H^2({\cal E}(H))=0$ and $\chi (%
{\cal E}(H))\geq \chi ({\cal E})$, then for all $n\geq 2$ the sheaf ${\cal F}%
(nH)$ is generated by global sections and its general section has degeneracy
locus of codimension $2$.
\end{lemma}

\TeXButton{Proof}{\proof}(i) By semicontinuity it is enough to exhibit one
such ${\cal F}$. We go by induction on $c_2$. If $c_2=c_2({\cal E})$, we may
take ${\cal F}={\cal E}$. If $c_2>c_2({\cal E})$, let ${\cal G}$ be a
prioritary sheaf of rank $r$ and Chern classes $c_1$ and $c_2-1$ with $H^1(%
{\cal G})=H^2({\cal G})=0$. By (\ref{RR}) we have $\chi ({\cal G})=\chi (%
{\cal F})+1>0$. So ${\cal G}$ must have a nonzero global section $s$. If $%
x\in S$ is a general point of $S$ and ${\cal G}\otimes k(x)\TeXButton{-->>}
{\twoheadrightarrow}k(x)$ a general one-dimensional quotient of the fiber of
${\cal G}$ at $x$, then the image of $s$ in $k(x)$ is nonzero. So if ${\cal F%
}$ is the kernel%
$$
0\rightarrow {\cal F}\rightarrow {\cal G}\rightarrow k(x)\rightarrow 0,
$$
then $h^0({\cal F})=h^0({\cal G})-1$ and $H^1({\cal F})=H^2({\cal F})=0$.

(ii) Under the added hypotheses the general ${\cal F}$ also satisfies $H^1(%
{\cal F}(H))=H^2({\cal F}(H))=0$. But $H^1({\cal F}(H))=H^2({\cal F})=0$
implies that ${\cal F}(nH)$ is generated by global sections for all $n\geq 2$
by the Castelnuovo-Mumford lemma. The general section of ${\cal F}(nH)$ will
drop rank in codimension $2$ by Bertini's theorem.

\begin{lemma}
\label{goodp2}Let ${\cal F}$ be a generic prioritary sheaf of rank $2$ and
Chern classes $c_1\geq -4$ and $c_2$ on ${\bf P}^2$. The two conditions (a) $%
H^1({\cal F})=H^2({\cal F})=0$ and (b) ${\cal F}(3)$ has a section with
degeneracy locus of codimension $2$ hold if and only if $\chi ({\cal F})\geq
0$.
\end{lemma}

\TeXButton{Proof}{\proof}If (a) and (b) hold, then clearly $\chi ({\cal F}%
)=h^0({\cal F})\geq 0$. Conversely, suppose ${\cal F}$ is generic prioritary
of rank $2$ with $c_1\geq -4$ and $\chi ({\cal F})\geq 0$. If $c_1=2a$ is
even, then let ${\cal E}={\cal O}_{{\bf P}^2}(a)^2$. Then $c_1=c_1({\cal E})$
and $c_2\geq \frac 14c_1^2=a^2=c_2({\cal E})$ by Proposition \ref{discp2}.
Moreover, $H^1({\cal E})=H^2({\cal E})=H^1({\cal E}(1))=H^2({\cal E}(1))=0$
and $\chi ({\cal E}(1))=\chi ({\cal E})+2a+4\geq \chi ({\cal E})$. Hence
conditions (a) and (b) follow from Lemma \ref{diminution}.

If $c_1=2a+1$ is odd, we may apply Lemma \ref{diminution} with ${\cal E}=%
{\cal O}_{{\bf P}^2}(a)\oplus {\cal O}_{{\bf P}^2}(a+1)$.\TeXButton{qed}
{\qed \medskip}

\begin{lemma}
\label{goodquad}Let ${\cal F}$ be a generic prioritary sheaf of rank $2$ and
Chern classes $c_1=(a_1,a_2)$ and $c_2$ on ${\bf P}^1\times {\bf P}^1$.
Suppose the $a_i\geq -2$.

(i)\quad If $a_2$ is even, then the two conditions (a) $H^1({\cal F})=H^2(%
{\cal F})=0$ and (b) ${\cal F}(2,2)$ has a section with degeneracy locus of
codimension $2$ hold if and only if $\chi ({\cal F})\geq 0$.

(ii)\quad If $a_2$ is odd, then (a) and (b) hold if and only if one has both
$\chi ({\cal F})\geq 0$ and either $c_2\geq \frac 12a_1(a_2-1)-1$ or $%
(a_1,a_2,c_2)=(a_1,-1,-a_1-2)$.
\end{lemma}

\TeXButton{Proof}{\proof}(i) We apply Lemma \ref{diminution} with ${\cal E}$
either ${\cal O}(\frac{a_1}2,\frac{a_2}2)^2$ or ${\cal O}(\frac{a_1-1}2,
\frac{a_2}2)\oplus {\cal O}(\frac{a_1+1}2,\frac{a_2}2)$.

(ii) Before beginning, recall that if $a_2=2b-1$ is odd, then by (\ref{Beil}%
) the general prioritary sheaf of the given rank and Chern classes is of the
form%
$$
0\rightarrow {\cal O}(a_1-p,b)\rightarrow {\cal F\rightarrow O}%
(p,b-1)\rightarrow 0
$$
with $p$ determined by $c_2=(a_1-p)(b-1)+pb=\frac 12a_1(a_2-1)+p$.

Now suppose that (a) and (b) hold. Then clearly $\chi ({\cal F})=h^0({\cal F}%
)\geq 0$. If $p\geq -1$, then $c_2\geq \frac 12a_1(a_2-1)-1$ as desired. If
on the other hand $p\leq -2$, then the sequence splits. So if (a) and (b)
hold, then $H^1({\cal O}(p,b-1))=0$ while ${\cal O}(p+2,b+1)$ is generated
by global sections. These are possible simultaneously only if $p=-2$ and $%
b=0 $.

Conversely, if $\chi ({\cal F})\geq 0$ and $c_2\geq \frac 12a_1(a_2-1)-1$,
then (a) and (b) hold by Lemma \ref{diminution} using ${\cal E}={\cal O}%
(a_1+1,b)\oplus {\cal O}(-1,b-1)$. If $(a_1,a_2,c_2)=(a_1,-1,-a_1-2)$, then
one may pick ${\cal F}={\cal O}(a_1+2,0)\oplus {\cal O}(-2,-1)$.
\TeXButton{qed}{\qed \medskip}

\section{Nonprioritary Components}

In this section we study nonprioritary components of ${\rm TF}_S(r,c_1,c_2)$
and of $W_N^0(E)$. According to Proposition \ref{unobst} on a ruled surface $%
\pi {:}~S\rightarrow C$ or on ${\bf P}^2$ with $f$ denoting the numerical
class either of a fiber of $\pi $ or of a line in ${\bf P}^2$, the general
member of any nonprioritary component of ${\rm TF}_S(r,c_1,c_2)$ a {\em %
nonprioritary generic extension of twisted ideal sheaves,} i.e.\ an
extension
\begin{equation}
\label{HN}0\rightarrow {\cal I}_{Z_1}(L_1)\rightarrow {\cal E}\rightarrow
{\cal I}_{Z_2}(L_2)\rightarrow 0
\end{equation}
such that the ${\cal O}_S(L_i)$ are generic line bundles having $L_1f>L_2f+1$
and the $Z_i$ are generic sets of $n_i$ points in $S$. In addition, the
proposition says that
\begin{equation}
\label{chi}\chi ({\cal I}_{Z_1}(L_1),{\cal I}_{Z_2}(L_2))=\chi ({\cal O}%
(L_2-L_1))-n_1-n_2\geq 0.
\end{equation}
Moreover, the extension is uniquely determined by ${\cal E}$ since it
defines the Harder-Narasimhan filtration of ${\cal E}$ with respect to a
suitable polarization of the surface.

The next two lemmas show that if $S$ is ${\bf P}^2$ or a semistable ruled
surface, then a generic extension of twisted ideal sheaves satisfying (\ref
{chi}) is the generic sheaf of an irreducible component of the stack of
torsion-free rank $2$ sheaves on $S$.

\begin{lemma}
\label{divisor}Suppose either that $\pi {:}~S\rightarrow C$ is a
birationally ruled surface or $S$ is ${\bf P}^2$. If a nonprioritary generic
extension of twisted ideal sheaves ${\cal E}$ as in (\ref{HN}) specializes
to another nonprioritary generic extension of twisted ideal sheaves $%
0\rightarrow {\cal I}_{Z_1^{\prime }}(L_1^{\prime })\rightarrow {\cal E}%
^{\prime }\rightarrow {\cal I}_{Z_2^{\prime }}(L_2^{\prime })\rightarrow 0$,
then

(i)\quad $\chi ({\cal I}_{Z_1^{\prime }}(L_1^{\prime }),{\cal I}%
_{Z_2^{\prime }}(L_2^{\prime }))<\chi ({\cal I}_{Z_1}(L_1),{\cal I}%
_{Z_2}(L_2))$ and

(ii)\quad there exists an effective divisor $\Gamma $ on $S$ such that $%
-\Gamma \cdot \Gamma >(L_1-L_2+K_S)\cdot \Gamma $.
\end{lemma}

\TeXButton{Proof}{\proof}(i) Let ${\rm FiltCoh}_S$ be the stack
parametrizing filtered coherent sheaves ${\cal F}_1\subset {\cal F}.$
Because the tangent space for automorphisms of ${\cal F}_1\subset {\cal F}$
(resp.\ the tangent space for deformations of ${\cal F}_1\subset {\cal F}$,
resp.\ the obstruction space for deformations of ${\cal F}_1\subset {\cal F}$%
) is ${\rm Ext}_{-}^i({\cal F},{\cal F})$ for $i=0$ (resp.$\ i=1$, resp.\ $%
i=2)$, one has%
$$
-\chi _{-}({\cal F},{\cal F})\leq \dim {}_{[{\cal F}_1\subset {\cal F}]}{\rm %
FiltCoh}_S\leq -\chi _{-}({\cal F},{\cal F})+\dim {\rm Ext}_{-}^2({\cal F},%
{\cal F})
$$
where $\chi _{-}({\cal F},{\cal F})=\sum (-1)^i\dim {\rm Ext}_{-}^i({\cal F},%
{\cal F})$. The forgetful functor ${\rm FiltCoh}_S\rightarrow {\rm Coh}_S$
defined by $[{\cal F}_1\subset {\cal F}]\mapsto [{\cal F}]$ induces maps on
infinitesimal automorphism, tangent, and obstruction spaces ${\rm Ext}_{-}^i(%
{\cal F},{\cal F})\rightarrow {\rm Ext}^i({\cal F},{\cal F})$. So if ${\rm %
Hom}_{+}({\cal F},{\cal F}):={\rm Hom}({\cal F}_1,{\cal F}/{\cal F}_1)=0$,
then the morphism ${\rm FiltCoh}_S\rightarrow {\rm Coh}_S$ is unramified at $%
[{\cal F}_1\subset {\cal F}]$, and ${\rm FiltCoh}_S$ can be viewed as more
or less a locally closed substack of ${\rm Coh}_S$ in a neighborhood of $[%
{\cal F}]$. In our case the subsheaf ${\cal F}_1={\cal I}_{Z_1}(L_1)$ is
unique, so ${\rm FiltCoh}_S$ is a locally closed substack of ${\rm Coh}_S$
in a neighborhood of $[{\cal F}]$.

Thus the dimension of the locally closed substack of torsion-free sheaves
numerically equivalent to ${\cal E}$ which admit a filtration with the
subsheaf numerically equivalent to ${\cal I}_{Z_1}(L_1)$ (resp.\ to ${\cal I}%
_{Z_1^{\prime }}(L_1^{\prime })$) and with ${\rm Ext}_{-}^2({\cal E},{\cal E}%
)=0$ is%
$$
-\chi _{-}({\cal E},{\cal E)=}-\chi ({\cal E},{\cal E})+\chi ({\cal I}%
_{Z_1}(L_1),{\cal I}_{Z_2}(L_2))
$$
(resp.\ $-\chi ({\cal E},{\cal E})+\chi ({\cal I}_{Z_1^{\prime
}}(L_1^{\prime }),{\cal I}_{Z_2^{\prime }}(L_2^{\prime }))$). If the former
substack contains the latter in its closure, its dimension must be larger.

(ii) As ${\cal E}$ specializes to ${\cal E}^{\prime }$, its subsheaf ${\cal I%
}_{Z_1}(L_1)$ specializes to a subsheaf of ${\cal E}^{\prime }$. Because
this subsheaf destabilizes ${\cal E}^{\prime }$, it must be contained in $%
{\cal I}_{Z_1^{\prime }}(L_1^{\prime })$. Hence ${\cal O}_S(L_1)$
specializes to a line bundle of the form ${\cal O}_S(L_1^{\prime }-\Gamma )$
with $\Gamma $ an effective divisor, and\ ${\cal O}_S(L_2)$ specializes to $%
{\cal O}_S(L_2^{\prime }+\Gamma )$.

Since $L_2^{\prime }-L_1^{\prime }\equiv L_2-L_1-2\Gamma $, the Riemann-Roch
formula leads to%
$$
\chi ({\cal O}(L_2^{\prime }-L_1^{\prime }))=\chi ({\cal O}(L_2-L_1))+\left(
2(L_1-L_2+\Gamma )+K_S\right) \cdot \Gamma \text{.}
$$
We also have
$$
n_1+n_2+\left( L_1\cdot L_2\right) =c_2({\cal E})=c_2({\cal E}^{\prime
})=n_1^{\prime }+n_2^{\prime }+\left( L_1^{\prime }\cdot L_2^{\prime
}\right) ,
$$
from which we see that%
$$
n_1^{\prime }+n_2^{\prime }=n_1+n_2+(L_1-L_2+\Gamma )\cdot \Gamma .
$$
Thus%
$$
\chi ({\cal I}_{Z_1^{\prime }}(L_1^{\prime }),{\cal I}_{Z_2^{\prime
}}(L_2^{\prime }))=\chi ({\cal I}_{Z_1}(L_1),{\cal I}_{Z_2}(L_2))+(L_1-L_2+%
\Gamma +K_S)\cdot \Gamma .
$$
Because of (i) this now implies the lemma.\TeXButton{qed}{\qed \medskip}

\begin{lemma}
\label{tilt}Suppose either that $\pi {:}~S\rightarrow C$ is a ruled surface
without curves of negative self-intersection or that $S$ is ${\bf P}^2$. If
a nonprioritary generic extension of twisted ideal sheaves ${\cal E}$ as in (%
\ref{HN}) specializes to another generic extension of twisted ideal sheaves $%
{\cal E^{\prime }}$, then $\chi ({\cal O}(L_2-L_1))\leq 0$.
\end{lemma}

\TeXButton{Proof}{\proof}Because ${\cal E}$ is not prioritary, the
restriction of ${\cal O}(L_2-L_1)$ to a general fiber of $\pi $ or a general
line of ${\bf P}^2$ is of negative degree. So $H^0({\cal O}(L_2-L_1))=0$.
Since $S$ contains no curves of negative self-intersection, Lemma \ref
{divisor}(ii) says that there is an effective divisor $\Gamma $ on $S$ such
that $(L_1-L_2+K_S)\cdot \Gamma <0$. Since $S$ contains no curves of
negative self-intersection, ${\cal O}(L_1-L_2+K_S)$ cannot be effective.
Thus $H^0({\cal O}(L_1-L_2+K_S))=0$ and by Serre duality $H^2({\cal O}%
(L_2-L_1))=0$. It follows that $\chi ({\cal O}(L_2-L_1))=-h^1({\cal O}%
(L_2-L_1))\leq 0$ as asserted. \TeXButton{qed}{\qed \medskip}

\begin{lemma}
\label{nonprior}Suppose either that $\pi {:}~S\rightarrow C$ be a ruled
surface without curves of negative self-intersection and $f\in {\rm NS}(S)$
is the class of a fiber of $\pi $, or that $S$ is ${\bf P}^2$ and $f$ is the
class of a line. Let $c_1\in {\rm NS}(S)$ and $c_2\in {\bf Z}$. Let $%
(D,n_1,n_2)\in {\rm NS}(S)\times {\bf Z}^2$. Then ${\rm TF}_S(2,c_1,c_2)$
has a unique component whose general member ${\cal E}$ is a nonprioritary
generic extension of twisted ideal sheaves
\begin{equation}
\label{again}0\rightarrow {\cal I}_{Z_1}(c_1-D)\rightarrow {\cal E}%
\rightarrow {\cal I}_{Z_2}(D)\rightarrow 0
\end{equation}
with $\deg (Z_i)=n_i$ if and only if $Df\leq -1+\frac 12c_1f$, and the $n_i$
are nonnegative and satisfy $n_1+n_2=c_2-D(D-c_1)\leq \chi ({\cal O}%
_S(2D-c_1))$. Such a component of ${\rm TF}_S(2,c_1,c_2)$ has dimension $%
-\chi ({\cal E},{\cal E})+\chi ({\cal O}_S(2D-c_1))-n_1-n_2$ and generic
embedding codimension $n_2+h^1({\cal O}_S(2D-c_1))$.
\end{lemma}

\TeXButton{Proof}{\proof}If ${\cal E}$ is a generic nonprioritary sheaf,
then its restriction to a general fiber $F$ of $\pi $ (resp.\ to a generic
line of ${\bf P}^2$) must be of the form ${\cal E}{\mid }_F\cong {\cal O}%
_F(a)\oplus {\cal O}_F(b)$ with $a\geq b+2$ by Lemma \ref{fiber} (resp. by
\cite{HL} Chap. I, Proposition 1.2). Hence the relative Harder-Narasimhan
filtration of ${\cal E}$ which was described before Lemma \ref{fiber}
(resp.~the Harder-Narasimhan filtration of ${\cal E}$ on ${\bf P}^2$) must
be of the form $0\subset {\cal I}_{Z_1}(c_1-D)\subset {\cal E}$ with ${\cal E%
}/{\cal I}_{Z_1}(c_1-D)\cong {\cal I}_{Z_2}(D)$ for some divisor $D$ on $S$
and some $0$-dimensional subschemes $Z_i\subset S$ such that $(c_1-D)f\geq
Df+2$, or $Df\leq -1+\frac 12c_1f$. Clearly one has $n_i:=\deg (Z_i)\geq 0$
and $c_2=D(c_1-D)+n_1+n_2$. Lemma \ref{fiber} shows that Proposition \ref
{unobst} is applicable to the filtered sheaf $0\subset {\cal I}%
_{Z_1}(c_1-D)\subset {\cal E}$. So $\chi ({\cal I}_{Z_1}(c_1-D),{\cal I}%
_{Z_2}(D))=\chi ({\cal O}_S(2D-c_1))-n_1-n_2\geq 0$. Thus to any
nonprioritary irreducible component of ${\rm TF}_S(2,c_1,c_2)$ there is an
associated triple $(D,n_1,n_2)$ satisfying the asserted numerical conditions.

Conversely suppose $(D,n_1,n_2)$ satisfy all the numerical conditions. Let $%
Z_i$ be a general set of $n_i$ points on $S$ and let ${\cal E}$ be a generic
extension as in (\ref{again}). Then ${\cal E}$ cannot be a specialization of
another nonprioritary generic extension $0\rightarrow {\cal I}_{Z_1^{\prime
}}(c_1-D^{\prime })\rightarrow {\cal E^{\prime }}\rightarrow {\cal I}%
_{Z_2^{\prime }}(D^{\prime })\rightarrow 0$ because in that case Lemma \ref
{divisor}(i) would imply%
$$
\chi ({\cal O}_S(2D^{\prime }-c_1))-n_1^{\prime }-n_2^{\prime }>\chi ({\cal O%
}_S(2D-c_1))-n_1-n_2\geq 0
$$
contradicting Lemma \ref{tilt}. Nor can ${\cal E}$ be a specialization of a
generic prioritary sheaf because it is the sheaf corresponding to a generic
point of a locally closed substack of ${\rm TF}_S(2,c_1,c_2)$ whose
dimension was calculated in the proof of Lemma \ref{divisor}(i) as $-\chi (%
{\cal E},{\cal E})+\chi ({\cal O}_S(2D-c_1))-n_1-n_2$. This is at least $%
-\chi ({\cal E},{\cal E})$, the dimension of the prioritary component. So $%
{\cal E}$ is the generic sheaf of an irreducible component of ${\rm TF}%
_S(2,c_1,c_2)$ of dimension $-\chi ({\cal E},{\cal E})+\chi ({\cal O}%
_S(2D-c_1))-n_1-n_2$.

The embedding codimension is the dimension of the cokernel of the map $%
\alpha $ between the tangent spaces of the stack of filtered sheaves and the
stack of unfiltered sheaves which is given by%
$$
{\rm Ext}_{-}^1({\cal E},{\cal E})\stackrel{\alpha }{\rightarrow }{\rm Ext}%
^1({\cal E},{\cal E})\rightarrow {\rm Ext}^1({\cal I}_{Z_1}(c_1-D),{\cal I}%
_{Z_2}(D))\rightarrow 0.
$$
But since ${\rm Hom}({\cal I}_{Z_1}(c_1-D),{\cal I}_{Z_2}(D))=0$, the
dimension of ${\rm cok}(\alpha )$ is the difference between the two numbers
\begin{eqnarray*}
& \dim {\rm Ext}^2({\cal I}_{Z_1}(c_1-D),{\cal I}_{Z_2}(D))=h^0({\cal I}%
_{Z_1}(c_1-2D+K_S))=\left[ h^2({\cal O}(2D-c_1))-n_1\right] _{+}, &  \\
& \chi ({\cal I}_{Z_1}(c_1-D),{\cal I}_{Z_2}(D))=\left[ h^2({\cal O}%
(2D-c_1))-n_1\right] -\left[ h^1({\cal O}(2D-c_1))+n_2\right] . &
\end{eqnarray*}
Because the $\chi $ is nonnegative, we see that $h^2({\cal O}%
(2D-c_1))-n_1\geq 0$, and that therefore the difference between the two
numbers is $n_2+h^1({\cal O}(2D-c_1))$. \TeXButton{qed}{\qed \medskip}

\begin{remark}
The components of ${\rm TF}_{{\bf P}^2}(2,c_1,c_2)$ containing locally free
sheaves were already classified by Str\o mme in \cite{St} Theorem 3.9, but
he made one minor error with the embedding codimensions. The prioritary
components of ${\rm TF}_{{\bf P}^2}(2,c_1,\frac 14c_1^2+1)$ are generically
smooth like all prioritary components. But they appear in Str\o mme's
classification in \cite{St} Theorem 3.9 as the component with $%
(d,c_1,c_2)=(0,0,1)$ which was said to be nonreduced with generic embedding
codimension $1$. The computation of the $h^i({\cal E}nd({\cal E}))$ in \cite
{St} Proposition 1.4 is wrong in that single case.
\end{remark}

We now consider what the classification of generic rank $2$ sheaves entails
for Brill-Noether loci. For the sake of simplicity, we will restrict
ourselves to those surfaces covered by Lemma \ref{tilt} which also have
vanishing irregularity, thus ${\bf P}^2$ and ${\bf P}^1\times {\bf P}^1$, so
that we do not need to analyze nongeneric line bundles which might have more
cohomology than the corresponding generic line bundles.

\begin{lemma}
\label{final}Let $S$ be ${\bf P}^1\times {\bf P}^1$ (resp.\ ${\bf P}^2$) and
let $f\in {\rm NS}(S)$ be the class of a fiber of ${\rm pr}_1$ (resp.\ a
line). Let ${\cal F}$ be a nonprioritary generic extension of twisted ideal
sheaves
\begin{equation}
\label{third}0\rightarrow {\cal I}_{Z_1}(c_1-D)\rightarrow {\cal F}%
\rightarrow {\cal I}_{Z_2}(D)\rightarrow 0
\end{equation}
with $c_1-K_S$ effective, $Df\leq -1+\frac 12c_1f$, and $\chi ({\cal O}%
_S(2D-c_1))\geq 0$. Let $n_i:=\deg (Z_i)$. Then the two conditions (a) $H^1(%
{\cal F})=H^2({\cal F})=0$ and (b) ${\cal F}(-K_S)$ has a section with
degeneracy locus of codimension $2$ hold if and only if the three conditions
hold: (i) $\chi ({\cal F})\geq 0$, (ii) $D-K_S$ is an effective and
irreducible divisor class, and (iii) $n_2\leq h^0({\cal O}_S(D))$.
\end{lemma}

\TeXButton{Proof}{\proof}We will prove the lemma for ${\bf P}^1\times {\bf P}%
^1$ only. The proof for ${\bf P}^2$ is similar and actually simpler.

Let $(a_1,a_2)$ be the bidegree of $c_1-D$ and $(b_1,b_2)$ the bidegree of $%
D $. We claim that the hypotheses of the lemma imply that $a_1\geq 0$ and $%
a_2\geq 0$. To see this first note that the effectiveness of $c_1-K_S$ is
equivalent to $a_1+b_1\geq -2$ and $a_2+b_2\geq -2$. The condition $Df\leq
-1+\frac 12c_1f$ is equivalent to $b_2\leq -1+\frac 12(a_2+b_2)$, or $%
a_2-b_2\geq 2$. Adding gives $a_2\geq 0$ as claimed. Also we have%
$$
0\leq \chi ({\cal O}_S(2D-c_1))=(b_1-a_1+1)(b_2-a_2+1).
$$
Since $b_2-a_2+1<0$, this gives $b_1-a_1+1\leq 0$. Thus $a_1>b_1$. Adding
this to $a_1+b_1\geq -2$ now gives $a_1\geq 0$ as claimed.

The fact that $c_1-D$ has bidegree $(a_1,a_2)$ with $a_1\geq 0$ and $a_2\geq
0$ implies that $H^i({\cal O}_S(c_1-D))=0$ for $i=1,2$ and that ${\cal O}%
_S(c_1-D+K_S)\TeXButton{ncong}{\ncong}{\cal O}_S$.

Now suppose that (a) and (b) hold. Then $\chi ({\cal F})=h^0({\cal F})\geq 0$%
, whence (i). To prove conditions (ii) and (iii), note first that $H^1({\cal %
F})$ and $H^2({\cal F})$ vanish because of (a) while $H^2({\cal I}%
_{Z_2}(c_1-D))\cong H^2({\cal O}_S(c_1-D))$ vanishes because of the previous
paragraph. So $H^1({\cal I}_{Z_2}(D))=H^2({\cal I}_{Z_2}(D))=0$ by (\ref
{third}). This implies that $H^1({\cal O}_S(D))=H^2({\cal O}_S(D))=0$ and $%
n_2\leq h^0({\cal O}_S(D))$. Thus we have (iii) plus $H^1({\cal O}%
_S(b_1,b_2))=H^2({\cal O}_S(b_1,b_2))=0$. These vanishings imply either that
both $b_i\geq -1$ and hence that the divisor $D-K_S$ of bidegree $%
(b_1+2,b_2+2)$ is very ample, or that $(b_1,b_2)$ is $(-1,d)$ or $(d,-1)$
with $d\leq -2$. But if $(b_1,b_2)$ had of one of these last two forms, and
if also $d\leq -3$, then ${\cal O}_S(D-K_S)={\cal O}_S(b_1+2,b_2+2)$ would
not have any global sections. Hence all sections of ${\cal F}(-K_S)$ would
lie in ${\cal I}_{Z_1}(c_1-D-K_S)$. But we have shown that the line bundle $%
{\cal O}_S(c_1-D-K_S)$ is always nontrivial. So all global sections of $%
{\cal F}(-K_S)$ would degenerate along a nontrivial curve, contradicting
(b). Hence the only possible cases where (a) and (b) hold with $D-K_S$ not
very ample are the cases where $D-K_S$ is of bidegree $(0,1)$ or $(1,0)$,
whence (ii). Thus (a) and (b) imply (i), (ii) and (iii).

Conversely, suppose (i), (ii) and (iii) hold for ${\cal F}$. We begin by
proving (a) in the special case where $n_1=0$. Condition (ii) implies that
either both $b_i\geq -1$ or one $b_i=-1$. Therefore $H^1({\cal O}_S(D))=H^2(%
{\cal O}_S(D))=H^2({\cal I}_{Z_2}(D))=0$. Because $Z_2$ consists of $n_2\leq
h^0({\cal O}_S(D))$ generic points of $S$ by condition (iii), we have $H^1(%
{\cal I}_{Z_2}(D))=0$ also. And we have already shown that $H^i({\cal O}%
_S(c_1-D))=0$ for $i=1,2$. It now follows by (\ref{third}) that $H^1({\cal F}%
)=H^2({\cal F})=0$. Thus (a) holds in the special case where $n_1=0$.

If $n_1>0$, then we may prove (a) by induction on $n_1$ using the same
method as in the proof of Lemma \ref{diminution}(i).

For (b) let $H$ be a divisor of bidegree $(1,1)$. Then%
$$
\chi ({\cal F}(H))=\chi ({\cal F})+(c_1+2H)H+2>\chi ({\cal F})\geq 0
$$
since $c_1+2H=c_1-K_S$ is effective. So (i) holds for ${\cal F}(H)$. Since
(ii) holds for ${\cal F}$, the divisor $D-K_S$ is base-point-free, so $%
D+H-K_S$ is very ample. Hence (ii) holds for ${\cal F}(H)$. And
$$
h^0({\cal O}_S(D+H))=h^0({\cal O}_S(D))+(D+H)H+1\geq h^0({\cal O}_S(D))\geq
n_2
$$
since according to (ii) $D+H$ is either effective or of bidegree $(0,-1)$ or
$(-1,0)$. So (iii) also holds for ${\cal F}(H)$. By what we have already
verified, the fact that (i), (ii) and (iii) all hold for ${\cal F}$ and $%
{\cal F}(H)$ implies that (a) also holds for ${\cal F}$ and ${\cal F}(H)$.
Hence $H^1({\cal F}(H))=H^2({\cal F})=0$. So ${\cal F}(2H)={\cal F}(-K_S)$
is generated by global sections by the Castelnuovo-Mumford lemma. Condition
(b) now follows from Bertini's theorem. \TeXButton{qed}{\qed \medskip}

\section{Proofs of the Theorems}

\paragraph{Proof of Theorems \ref{ruled} and \ref{P2}.}

Theorem \ref{ruled} follows from the classification of the prioritary
components of ${\rm TF}_S(2,c_1,c_2)$ in Lemma \ref{discquad} and the
classification of the nonprioritary components of ${\rm TF}_S(2,c_1,c_2)$ in
Lemma \ref{nonprior}. Note that the expression $\chi ({\cal O}%
_S(-c_1))+D(2D-2c_1-K_S)$ appearing the Theorem \ref{ruled} is equal to the
expression $\chi ({\cal O}_S(2D-c_1))$ appearing in Lemma \ref{nonprior} by
a simple application of the Riemann-Roch formula for a line bundle on $S$.

Theorem \ref{P2} follows from Lemmas \ref{discp2} and \ref{nonprior} in the
same manner.

\paragraph{Proof of Theorem \ref{BN}.}

According to the argument given before the statement of Theorem \ref{BN}
there is a correspondence between irreducible components of $W_N^0(E)$
correspond to the irreducible components of ${\rm TF}_S(2,E-K_S,N)$ whose
general member ${\cal E}$ satisfies $H^1({\cal E}(K_S))=H^2({\cal E}(K_S))=0$
and has a section with zero locus of codimension $2$. These irreducible
components of ${\rm TF}_S(2,E-K_S,N)$ may either be nonprioritary or
prioritary. According to Theorems \ref{ruled} and \ref{P2} the nonprioritary
components of ${\rm TF}_S(2,E-K_S,N)$ correspond to pairs $(D,n)\in {\rm NS}%
(S)\times {\bf Z}$ such that $Df\leq -1+\frac 12(E-K_S)f$ and $0\leq n\leq
N+D(D-E+K_S)$ and $N\leq \chi ({\cal O}_S(-E+K_S))+D(D-E)=\chi ({\cal O}%
_S(E))-D(E-D)$. According to Lemma \ref{final} the general member of such an
irreducible component has $H^1({\cal E}(K_S))=H^2({\cal E}(K_S))=0$ and a
section with zero locus of codimension $2$ if and only if (i) $\chi ({\cal E}%
(K_S))\geq 0$, (ii) $D$ is an effective and irreducible divisor class, and
(iii) $n_2=N+D(D-E+K_S)-n\leq \chi ({\cal O}_S(D+K_S))$. Since $\chi ({\cal E%
}(K_S))=\chi ({\cal I}_X(E))+1=\chi ({\cal O}_S(E))-N+1>0$, we see that the
irreducible components of $W_N^0(E)$ with nonprioritary ${\cal E}$ are
precisely the components described in part (i) of Theorem \ref{BN}.
Moreover, the geometry of $X$ can be recovered from ${\cal E}$, $D$ and $n$
via the diagram%
$$
\begin{array}{ccccccc}
&  &  & 0 &  & 0 &  \\
&  &  & \downarrow &  & \downarrow &  \\
&  &  & {\cal O}_S & = & {\cal O}_S &  \\
&  &  & \downarrow &  & \downarrow &  \\
0\rightarrow & {\cal I}_{Z_1}(E-K_S-D) & \rightarrow & {\cal E} &
\rightarrow & {\cal I}_{Z_2}(D) & \rightarrow 0 \\
& \parallel &  & \downarrow &  & \downarrow &  \\
0\rightarrow & {\cal I}_{Z_1}(E-K_S-D) & \rightarrow & {\cal I}_X(E-K_S) &
\rightarrow & {\cal K} & \rightarrow 0 \\
&  &  & \downarrow &  & \downarrow &  \\
&  &  & 0 &  & 0 &
\end{array}
$$
The bottom row must be a twist of the residual exact sequence for ${\cal I}%
_X(E-K_S)$ with respect to a curve $C\in |D|$. So ${\cal K}={\cal I}_{X\cap
C/C}(E-K_S)$. Thus $X=Z_1\cup (X\cap C)$ with $Z_1$ a generic set of $n$
points of $S$ and $X\cap C$ a set of $N-n$ points on $C$. Thus the
irreducible components of $W_N^0(E)$ such that ${\cal E}$ is nonprioritary
are exactly those described in part (i) of Theorem \ref{BN}.

In addition ${\rm TF}_S(2,E-K_S,N)$ may have a unique prioritary component.
For ${\bf P}^2$ this component exists if and only if $N\geq \frac
14(e+4)(e+2)$ by Theorem \ref{P2}. Its general member ${\cal E}$ always
satisfies $H^1({\cal E}(K_S))=H^2({\cal E}(K_S))=0$ and has a section with
zero locus of codimension $2$ according to Lemma \ref{goodp2} because $\chi (%
{\cal E}(K_S))=\chi ({\cal O}_S(E))-N+1>0$. If $S$ is ${\bf P}^1\times {\bf P%
}^1$ and $e_2$ is even, the prioritary component exists if and only if $%
N\geq \frac 12(e_1+2)(e_2+2)$ according to Theorem \ref{ruled}, and its
general member always satisfies $H^1({\cal E}(K_S))=H^2({\cal E}(K_S))=0$
and has a section with zero locus of codimension $2$ according to Lemma \ref
{goodquad} because $\chi ({\cal E}(K_S))>0$. If $S$ is ${\bf P}^1\times {\bf %
P}^1$ and $e_2$ is odd, then the prioritary component exists for all $N$
according to Theorem \ref{ruled}. But according to Lemma \ref{goodquad} its
general member ${\cal E}$ satisfies $H^1({\cal E}(K_S))=H^2({\cal E}(K_S))=0$
and has a section with zero locus of codimension $2$ only if either $N\geq
\frac 12(e_1+2)(e_2+1)+1$ or $(e_1,e_2,N)=(e_1,1,e_1+2)$. This gives all the
components described in part (ii) of Theorem \ref{BN}.

The dimension of a component of $W_N^0(E)$ is the dimension of the
corresponding component of ${\rm TF}_S(2,E-K_S,N)$ plus $h^0({\cal E})-1$
(for the choice of a section of ${\cal E}$ modulo $k^{\times }$) plus $1$
(to cancel the negative contribution of $\dim {\rm Aut}({\cal I}_X(E))$ in
the stack computations). Hence the prioritary components have dimension $%
-\chi ({\cal E},{\cal E})+\chi ({\cal E})$ which a straightforward
Riemann-Roch computation shows is $2N-(\chi ({\cal O}_S(E)-N+1)$. Since $%
\dim {\rm Hilb}^N(S)=2N$, this is the asserted codimension $(\chi ({\cal O}%
_S(E)-N+1)$. The nonprioritary components of ${\rm TF}_S(2,E-K_S,N)$ have
dimensions greater by $\chi ({\cal O}_S(E))+D(D-E)-N$. So the nonprioritary
components of $W_N^0(E)$ have codimensions $D(E-D)+1$. \TeXButton{qed}
{\qed \medskip}

\end{document}